\documentstyle[12pt,leqno]{article}
\title{On the Construction of Isospectral Vectorial \\
        Sturm-Liouville Differential Equations 
        }
\author{ Hua-Huai Chern }
\date{ March 25, 1998 }
\begin{document}
\baselineskip=20pt
\maketitle
\vskip -0.5cm
\centerline{\footnotesize Department of Mathmatics}
\centerline{\footnotesize National Chung Cheng University}
\centerline{\footnotesize Minghsiug, Chiayi 621}
\centerline{\footnotesize hhchern@math.ccu.edu.tw }

\vskip 0.5cm
\par
\vskip 2.5cm
\abstract
We extend the idea of Jodeit and Levitan \cite{JL} for constructing
isospectral problems of the classical scalar
Sturm-Liouville differential equations to the vectorial Sturm-Liouville differential equations.
Some interesting relations are found.
\vskip 0.5cm
\noindent {\bf Keywords and phrases.} {\rm vectorial Sturm-Liouville differential equation,
matrix differential equations, isospectral problem } \par
\vskip 0.5cm
\noindent {\bf AMS(MOS) subject classifications. } 34A30, 34B25. \par
\vskip 1cm
\noindent Abbreviated title: Isospectral Problems
\newpage

\section{ Introduction }

The purpose of this article is to study the problem of constructing
$N$-dimensional, $N \ge 2$ , vectorial Sturm-Liouville differential equations
subject to certain boundary conditions,
which are isospectral to a given one.

Those $N$-dimensional vectorial Sturm-Liouville eigenvalue problems which are
considered in this paper, are of the following form:
\begin{equation} \label{a1.1}
      -\phi ''(x) +P(x) \phi (x) = \lambda \phi (x), \quad
       B\phi '(0) + A \phi (0) = {\cal B} \phi ' (\pi )
          +{\cal A} \phi (\pi ) = {\bf 0 } ,
\end{equation}
where $0\le x \le \pi $ , $\phi (x) $ is an ${\bf R}^N$-valued function,
$P(x)$ is a continuous $N
\times N $ symmetric matrix-valued function, $A, B, {\cal A, B} $ are $N
\times N $ matrices which satisfy the following conditions:
\begin{equation} \label{a1.2}
BA^* = AB^* , {\cal BA^* = AB^* }, \quad \mbox{rank}[ A, B]=\mbox{rank}[ {\cal
A, B } ]= N ,
\end{equation}
where $A^*$ is the transpose matrix of $A$, $[A,B]$ denotes the $N \times 2N $
matrix whose first $N\times N$ block
is $A$, and the second $N \times N$ block is $B$. We shall use the tuple $(P, A,
B, {\cal A, B })$ to denote the eigenvalue problem (\ref{a1.1}). Note that the
conditions in (\ref{a1.2}) ensure the problem (\ref{a1.1}) a selfadjoint eigenvalue
problem, and its eigenvalues can be determined by the variational principle.
Counting multiplicities of the eigenvalues, we arrange the eigenvalues of (\ref{a1.1})
in an ascending sequence
\begin{equation} \label{a1.3}
      \mu _0 \le \mu _1 \le \mu _2 \le \cdots .
\end{equation}
This sequence shall be denote by $\Sigma (P, A, B, {\cal A, B}) $, and is called
the {\sl sequence of} \newline  {\sl eigenvalues} of (\ref{a1.1}).
Note that the multiplicity
of each eigenvalue of (\ref{a1.1}) is at most $N$. For convenience, we shall use
$\sigma (P,A,B,{\cal A,B})$ to denote the {\sl set of eigenvalues} of
(\ref{a1.1}), and arrange its elements in ascending order as
\[
      \lambda _0 < \lambda _1 < \lambda _2 < \cdots  ,
\]
and use $m_k$ denote the multiplicity of $\lambda _k $ in the sequence (\ref{a1.3}).
Given two $N$-dimensional selfadjoint vectorial Sturm-Liouville eigenvalue
problem $(P, A, B, {\cal A,} $ \newline ${\cal B}) $ and $(\tilde{P} ,\tilde{A} ,\tilde{B} ,\tilde{\cal A} , \tilde
{\cal B} )$ over $ 0 \le x \le \pi $.
If $\Sigma (P, A, B, {\cal A, B})=
\Sigma (\tilde{P} ,\tilde{A} ,\tilde{B} ,\tilde{\cal A} , \tilde{\cal B} )$, we
call these two eigenvalue problems {\sl isospectral problems}, or
$(P, A, B, {\cal A, B}) $ is {\sl isospectral} to
$(\tilde{P} ,\tilde{A} ,\tilde{B} ,\tilde{\cal A} , \tilde{\cal B} )$. For scalar
Sturm-Liouville equations, i.e., (\ref{a1.1}) with $N=1$, isospectral problems have
been studied by many mathematicians, notably, G. Borg \cite{B}, I. M. Gel'fand
, B. M. Levitan and their associates, and the structure of the set of isospectral
problems for scalar Sturm-Liouville equations is well-presented in the book of
J. P\"{o}schel and E. Trubowitz \cite{PT}, and in the works of E. Trubowitz. But
for vectorial Strum-Liouville equations, i.e., (\ref{a1.1}) with $N \ge 2$ ,
methods for constructing isospectral problems and the structure of isospectral problems are not
well-understood. Motivated by a recent work of Jodeit and Levitan \cite{JL},
in this paper we present a method for constructing an $N$-dimensional vectorial
Sturm-Liouville eigenvalue problem
$(\tilde{P} ,\tilde{A} ,\tilde{B} ,\tilde{\cal A} , \tilde{\cal B} )$
over $ 0 \le x \le \pi $, which is isospectral to the given $N$-dimensional
vectorial Sturm-Liouville eigenvalue problem $(P, A, B, {\cal A, B}) $
over $ 0 \le x \le \pi $. For simplicity we shall assume the completeness
of the eigenfunctions of the the given eigenvalue problem $(P, A, B, {\cal A, B}) $
subject to the given boundary conditions shown in (\ref{a1.1}).

This paper is organized as follows. In section 2 we present some preliminary
results related to vectorial Sturm-Liouville equations, and a matrix wave
equation (see (\ref{Th3-1})) which is constructed from the given eigenvalue problem
$(P, A, B, {\cal A, B}) $.
In section 3, using matrix wave equation introduced in section 2, we
construct eigenvalue problems which are isospectral to $(P, A, B, {\cal A, B}) $.
In section 4 we present two examples using our construction method. As shown in one
of our examples, even for a simple case such as $(P,I,0,I,0)$, where $P(x)$ is a
constant two by two diagonal matrix, $ I$ is the two by two identity matrix,
an isospectral problem $(Q,I,0,I,0)$ can be found where $Q(x)$ is not
simultaneously diagonalizable. The isospectral problem for vectorial Sturm-Liouville
equations is much more complicated than its scalar counterparts.

\section{ Preliminary }

To study the eigenvalue problem (\ref{a1.1}) we introduce the following matrix
differential equation
\begin{equation} \label{a2.1}
- Y'' +P(x)Y= \lambda Y, \mbox{} \hskip 1cm Y(0)=B^* ,
                         \mbox{} \hskip 1cm Y'(0)=-A^*.
\end{equation}
Let $Y(x; \lambda) $ denote the $N \times N$ matrix-valued solution of the
initial value problem (\ref{a2.1}). We have ( see \cite{CS} ),
\[
    Y(x; \lambda )= {\cal C} (x; \mu ) + \int _0^x \tilde{K} (x,t) {\cal C}
                                          (t; \mu ) dt \]
where $\mu ^2 = \lambda $ , $ {\cal C } (x; \mu ) = \cos (\mu x) B^* -
\mu^{-1} \sin (\mu x) A^* $, and  $\tilde{K}(x,t)$ is as that described
 in [2, {\bf Lemma 2.1}]. Define the following matrix-valued function
\begin{equation} \label{a2.2}
   W(\lambda ) = {\cal B} Y' (\pi ; \lambda )
                         +{\cal A} Y (\pi ;  \lambda ).
\end{equation}
Then $\lambda_* \in \sigma (P, A, B, {\cal A, B})$ if and only if $W (\lambda _*)$
is a singular matrix. It follows from the variational principle that the set of
the zeros of the equation
\begin{equation} \label{a2.3}
        \det W (\lambda ) = 0
\end{equation}
is bounded below. Denote the distint zeros of (\ref{a2.3}) in stricting ascending
order as
\begin{equation} \label{a2.4}
      \lambda _0 < \lambda _1 < \lambda _2 < \cdots .
\end{equation}
Then the multiplicity $m_k$ of $\lambda_k$ in the sequence of eigenvalues
$\Sigma (P, A, B, {\cal A, B})$ of (\ref{a1.1}) is equal to the dimemsion of
the null space $\mbox{Null}(W( \lambda_k ))$ of $W( \lambda _k )$. If ${\bf v}$ is a
nonzero element in $\mbox{Null}(W( \lambda_k ))$, then the vector-valued function
\[
         z(x) = Y (x; \lambda_k){\bf v}
\]
is an eigenfunction of (\ref{a1.1}) corresponding to the eigenvalue $ \lambda_k$.
In the following, for ${\bf v}$ and ${\bf w}$ in ${\bf R}^N$,
the notation $\langle {\bf v,w} \rangle $ is used to denote the inner
product of two elements ${\bf v}$ and ${\bf w}$ in ${\bf R}^N$.
We shall need the following result.
\newtheorem{c1}{Lemma}[section] \begin{c1} \label{le1}
For each $k \ge 0$, in the null space Null$(W( \lambda_k ))$ there are
exactly $m_k$ linearly independent constant vectors $\theta _l (k) $, $1\le
l\le m_k $, such that the vector-valued functions
\[   Y(x;\lambda_k)\theta _l (k) ,\mbox{} \hskip 1cm 1\le l \le m_k , \]
are mutually orthogonal, i.e.,
\[
   \int _0^\pi \langle Y(x; \lambda _k ) \theta _i (k) , Y(x; \lambda _k )
    \theta _j (k) \rangle dx = 0 ,  \hskip 0.5cm \mbox{if} \hskip 0.5cm i \neq
j. \]
\end{c1}
\noindent {\bf Proof.} Let $ v_1 , \ldots , v_{m_k} $ be a basis of
Null$(W( \lambda _k ))$ and $V_k =[ v_1 ,\ldots , v_{m_k} ]$. Then
\[
      V= \int_0^\pi V_k^* Y^* (x; \lambda _k)Y(x; \lambda _k )V_k dx
\]
is an $m_k \times m_k $ positive definite matrix. There exists an $m_k \times
m_k $ orthogonal matrix $U$ which diagonalizes $V$, i.e., $U^* V U $ is a
diagonal matrix. Let $\theta _l (k) , 1\le l \le m_k $, denote the column
vectors of $VU$. Then $ \theta _l (k) $ fulfills our requirement. $\Box $

According to {\bf Lemma \ref{le1} }, we define the following vector-valued
function
\begin{equation}
   \phi _l (x; \lambda , \lambda _k )= Y(x; \lambda ) \theta _l (k) ,
   \mbox{} \hskip 1.2cm  1 \le l \le m_k.
\end{equation}
Then, the functions
\[
   \phi _l (x; \lambda _k , \lambda _k )= Y(x; \lambda_k ) \theta _l
(k) ,    \mbox{} \hskip 1.2cm  1 \le l \le m_k,
\]
form an orthogonal basis of the eigenspace corresponding to the eigenvalue
$\lambda _k$ . From now on the eigenvalue problem $ (P,A,B,{\cal A,B})$
shall be fixed. And, as it was mentioned in the introduction, for simplicity,
we shall assume the completeness of the system of eigenfunctions
$\{ \phi _l (x; \lambda_k , \lambda_k ) : 1\le l \le m_k , k =0,1,2,\ldots \} $
subject to the given boundary conditions.

Next we extend the idea used by Jodeit and Levitan in \cite{JL} to construct the
kernel function for a related integral equation. We shall view an ${\bf R}^N$-vector
${\bf v}$ as an $N \times 1$ matrix.
 Choose $c_k^i \in R $, $1\le i \le m_k,
k=0,1,\ldots $, which convergs so rapidly to zero that the matrix-valued function ${\cal F}$ ,
defined by the following uniformly convergent series,
\begin{equation} \label{F}
     {\cal F}(x,y)= \sum_{k=0}^{\infty} \sum_{i=1}^{m_k}
                 c_k^i \phi _i (x;\lambda _k , \lambda _k )
                       \phi _i^* (y;\lambda _k , \lambda _k ),
\end{equation}
is continuous and has continuous first and second order derivatives. Then we
construct the following integral equation :
\begin{equation} \label{IE}
    K(x,y) + {\cal F}(x,y) + \int _0^x K(x,t) {\cal F}(t,y) dt =0
                    , \mbox{} \hskip 0.5cm 0\le y<x \le \pi.
\end{equation}

We note that if we choose the sequence $( c_k^i :1 \le i \le m_k, k=0,1,2, \ldots ) $
such that $c_k^i=0$ for all $ k \ge k_\circ $, $ i=1,2,\ldots ,m_k $, where $k_\circ $
is a fixed index, then the series in the right hand side of (\ref{F}) is a
finite series, and the equation (\ref{IE}) makes sense. This choice of $(c_k^i )$ shall
be used in section 4 to construct some concrete examples.

The existence of the solution of (\ref{IE}) can be easily proven by using
iteration method. On the other hand, when we choose suitably those real numbers
$c_k^i$,$1\le i \le m_k$ , $k \ge 0 $, we may prove the following uniqueness
theorem. \par

\newtheorem{c2}[c1]{Theorem} \begin{c2} \label{Th2.2}
Suppose that the sequence $( c_k^i )$ is chosen so that the series in (\ref{F})
is uniformly convergent and has continuous first and second order derivatives,
and
\begin{equation} \label{Th2-c}
      1+c_k^i || \phi _i (\cdot ; \lambda _k , \lambda _k )||^2 > 0 ,
       \mbox{} \hskip 0.5cm  1\le i \le m_k , \hskip 0.5cm \forall k \ge 0 .
\end{equation}
Then (\ref{IE}) has a unique solution for every $x$ , $ 0<x \le \pi $.
\end{c2}

\noindent {\bf Proof. } The method for proving this theorem is similar to the
one used for treating the scalar case in [ 3, {\bf Theorem 1.1}] .
 It suffices to show that the only solution for
the integral equation
\begin{equation} \label{IE2}
     \Delta (x,y) + \int _0^x \Delta (x,t) {\cal F}(t,y )dt =  0
\end{equation}
is $\Delta (x,y ) \equiv  0 $, where $\Delta (x,y )$ is the difference of two
solution of (\ref{IE}). Denote $\phi_{i,k} (x) = \phi_i (x; \lambda_k, \lambda_k)$
for convenience.
Owing to the assumption about the completeness of eigenfunctions of
$(P, A, B, {\cal A,B})$, and the
orthogonality ({\bf Lemma \ref{le1}}) of the eigenfunctions $\phi_{i,k} (x)$,
we have
\begin{equation}  \label{Th2.2-1}
   \Delta^* (x,y) = \sum_{k=0}^\infty \sum_{i=1}^{m_k}
                    \frac {\phi_{i,k} (y)}{||\phi_{i,k}||}
                    ( \int_0^x \frac {\phi_{i,k}^* (t)}{||\phi_{i,k}||}
                       \Delta^* (x,t) dt ) .
\end{equation}
By (\ref{IE2}), we have
\[
     \Delta^* (x,y) + \int _0^x {\cal F}^* (t,y )\Delta^* (x,t) dt =  0,
\]
\begin{equation} \label{Th2.2-2}
     \Delta (x,y)\Delta^* (x,y) + \int _0^x \Delta (x,y){\cal F}^* (t,y )\Delta^* (x,t) dt =  0.
\end{equation}
Integrating (\ref{Th2.2-2}) with respect to $y$-variable from $0$ to $x$, using (\ref{F})
and (\ref{Th2.2-1}), we obtain
\begin{equation} \label{Th2.2-3}
   \sum_{k=0}^\infty \sum_{i=0}^{m_k}
    \frac 1{||\phi_{i,k}||}[ 1+ c_k^i ||\phi_{i,k}||^2] \Omega_{i,k}(x)\Omega_{i,k}^*(x)=0,
\end{equation}
where
\begin{equation} \label{Th2.2-4}
\Omega_{i,k}(x)= \int_0^x \Delta (x,t) \phi_{i,k} (t) dt.
\end{equation}
Since $ \Omega_{i,k}(x)\Omega_{i,k}^*(x)$ is nonnegative definite,
(\ref{Th2.2-3}) and (\ref{Th2-c}) imply that
\[
  \Omega_{i,k}(x)\Omega_{i,k}^*(x)=0 ,
\]
and hence $ \Omega_{i,k}(x)=0 $, and by (\ref{Th2.2-4}),
\begin{equation} \label{Th2.2-5}
   \int_0^x \Delta (x,t) \phi_{i,k} (t) dt=0
\end{equation}
for $ k=0,1,2, \ldots, i= 1,2, \ldots ,m_k $. Then by the completeness of
eigenfunctions of $(P, A, B, {\cal A,B})$, (\ref{Th2.2-5}) implies
$\Delta (x,y) =0 $. This completes the proof. $\Box $

\vskip 0.5cm
Now we face the question :
`` Does the matrix-valued function $K(x,y)$ determined by the above theorem
also
satisfy some wave equation with which we are familiar as in the scalar case ? ''
The answer is affirmative, as shown below. \par

\newtheorem{c3}[c1]{Theorem} \begin{c3} \label{Th3}
Assumption as {\bf Theorem \ref{Th2.2}}.
The solution $K(x,y)$ of (\ref{IE}) satisfies the following partial
differential equation
\begin{equation} \label{Th3-1}
     \frac {\partial ^2}{\partial x^2} K -Q(x) K =
             \frac {\partial ^2}{\partial y^2} K -KP(y) ,
\end{equation}
where $ Q(x)= P(x)+2 d/dx K(x,x) $, and it also satisfies the following
conditions:
\[   K(x,y) =0 , \mbox{} \hskip 1cm \mbox{ if} \hskip 1cm y > x , \]
\begin{equation} \label{Th3-2}
            K(x,0)A^* + \frac \partial{\partial y}K|_{y=0}B^* = 0 ,
\end{equation}
\begin{equation} \label{Th3-3}
        K(x,x)= \frac 12 \int _0^x [ Q(t) -P(t)]dt -{\cal F}(0,0) ,
\end{equation}
where
\begin{equation} \label{Th3-4}
      {\cal F}(0,0) = B^* ( \sum_{k=0}^\infty \sum_{i=1}^{m_k}
                          c_k^i \theta _i(k) \theta _i^*(k) )B.
\end{equation}
\end{c3}

\noindent {\bf Proof. } Denote
\[
       {\cal J} (x,y)=
    K(x,y) + {\cal F}(x,y) + \int _0^x K(x,t) {\cal F}(t,y) dt .
\]
Then by (\ref{IE}), ${\cal J}=0$, hence ${\cal J}_{xx}={\cal J}_{yy}=0$. On the other hand, as
\begin{eqnarray*}
  {\cal J}_{xx} &=& \mbox{}
         \frac {\partial ^2}{\partial x ^2}K
         +[ P(x)
                    + ( \frac d{dx} K(x,x) + \frac {\partial }{\partial x}
                       K(x,x))]{\cal F} (x,y) \\
         & &\mbox{} - \sum_{k=0}^{\infty} \sum_{i=1}^{m_k} \lambda _k
                       c_k^i \phi _i (x;\lambda _k , \lambda _k )
                       \phi _i^* (x;\lambda _k , \lambda _k )    \\
         & &\mbox{} +K(x,x) \frac {\partial }{\partial x} {\cal F}(x,y)
                    +\int _0^x \frac {\partial ^2}{\partial x^2}K(x,t)
                      {\cal F}(t,y)dt,
\end{eqnarray*}
\begin{eqnarray*}
           {\cal J}_{yy}&=&\frac {\partial ^2 }{\partial y^2} K
                           +({\cal F}(x,y)+\int _0^x K(x,t){\cal F}(t,y)dt
                             )P(y) \\
                        & &\mbox{}+
                             \sum_{k=0}^{\infty} \sum_{i=1}^{m_k} \lambda _k
                              c_k^i \phi _i (x;\lambda _k , \lambda _k )
                              \phi _i^* (x;\lambda _k , \lambda _k )    \\
                        & &\mbox{}-\int _0^x K(x,t)[
                             \sum_{k=0}^{\infty} \sum_{i=1}^{m_k} \lambda _k
                              c_k^i \phi _i (t;\lambda _k , \lambda _k )
                              \phi _i^* (y;\lambda _k , \lambda _k )]dt,
\end{eqnarray*}
we have
\begin{eqnarray*}
    0 &=&\mbox{} {\cal J}_{xx}-{\cal J}_{yy} +{\cal J}P(y)  \\
      &=&\mbox{} \frac {\partial ^2}{\partial x^2} K-\frac {\partial ^2}{
                 \partial y^2}K +[ P(x)+ 2\frac d{dx}K(x,x)- \frac {\partial}{
                  \partial y} K|_{y=x} ]{\cal F}(x,y) \\
      & &\mbox{}+K(x,x) \frac {\partial }{\partial x} {\cal F}
                + \int _0^x \frac {\partial ^2}{\partial x^2} K(x,t){\cal F}
                  (t,y) dt \\
      & &\mbox{} +
              \int _0^x K(x,t)[
                            \sum_{k=0}^{\infty} \sum_{i=1}^{m_k} c_k^i
                            ( \lambda _k \phi _i (t;\lambda _k , \lambda _k ))
                             \phi _i^* (y;\lambda _k , \lambda _k )]dt.
\end{eqnarray*}
Replacing $\lambda_k \phi _i (t;\lambda _k , \lambda _k ) $ by $ -\phi_i ''
(t;\lambda _k ,\lambda _k ) +P(t) \phi _i (t; \lambda _k , \lambda _k ) $ in
the last integral and using integration by parts twice, we obtain
\begin{eqnarray*}
   0 &=&\mbox{} {\cal J}_{xx} - {\cal J}_{yy}+{\cal J}P(y) \\
     &=&\mbox{} \frac {\partial ^2}{\partial x^2} K
               -\frac {\partial ^2}{\partial y^2} K+(P(x)+2 \frac d{dx} )
                {\cal F}(x,y) \\
     & &\mbox{} +[ K(x,0) \frac {\partial}{\partial x}{\cal F}(x,y)|_{x=0}
                   -\frac {\partial }{\partial y} K(x,y)|_{y=0}{\cal F}
                   (0,y) ] \\
     & &\mbox{} + \int _0^x [ \frac {\partial ^2}{\partial x^2} K(x,t)-
                              \frac {\partial ^2}{\partial t^2} K(x,t)
                              +K(x,t)P(t)]{\cal F}(t,y) dt .
\end{eqnarray*}
Finally, let $ Q(x)=P(x)+2d/{dx} K(x,x) $, we have
\begin{eqnarray}
  0 &=&\mbox{} {\cal J}_{xx} -{\cal J}_{yy} +{\cal J}P(y)-Q(x){\cal J}
        \nonumber \\
    &=&\mbox{} +[ K(x,0) \frac {\partial}{\partial x}{\cal F}(x,y)|_{x=0}
                  -\frac {\partial }{\partial y} K(x,y)|_{y=0}{\cal F}
                  (0,y) ] \nonumber \\
    & & \mbox{}+[ \frac {\partial ^2}{\partial x^2} K
                  -\frac {\partial ^2}{\partial y^2} K -Q(x)K+KP(y)]
 \label{Th3-5}
\\     & &\mbox{}+\int _0^x
                [ \frac {\partial ^2}{\partial x^2} K
                  -\frac {\partial ^2}{\partial t^2} K -Q(x)K+KP(t)]
                {\cal F}(t,y)dt , \nonumber
\end{eqnarray}
where the function
\[
      K(x,0) \frac {\partial}{\partial x}{\cal F}(x,y)|_{x=0}
      -\frac {\partial }{\partial y} K(x,y)|_{y=0}{\cal F}
      (0,y)  \]
vanishes if and only if (\ref{Th3-2}) holds. Then (\ref{Th3-5}) becomes
an integral equation of the same type as (\ref{IE2}). Hence we have
(\ref{Th3-1}). \par

   By  (\ref{F}) and (\ref{IE}), we see that
\[
         K(0,0)=-{\cal F}(0,0) =B^*( \sum_{k=0}^\infty \sum_{i=1}^{m_k}
                                   c_k^i \theta _i(k) \theta _i^*(k) )B .
\]
(\ref{Th3-3}) is a consequence of the fundamental theorem of calculus
from the definition of $Q(x)$ given above.
$\Box $
\par

\newtheorem{c4}[c1]{Theorem} \begin{c4} \label{Th4}
If $K(x,y)$ is determined by {\bf Theorem \ref{Th3} }, then, for every
complex $\lambda $, $k\ge 0$ and $ 1\le l \le m_k $, the
vector-valued function $\psi _l (x; \lambda , \lambda_k )$ defined by
\begin{equation} \label{Th2.4-1}
   \psi _l (x; \lambda , \lambda_k )= \phi _l (x; \lambda , \lambda_k )
          + \int _0^x K(x,t) \phi _l (t; \lambda , \lambda_k ) dt
\end{equation}
is a solution of the vectorial differential equation
\begin{equation} \label{Th2.4-2}
      -\psi '' + Q(x) \psi = \lambda \psi , \mbox{} \hskip 1cm 0\le x\le \pi ,
\end{equation}
where $Q(x)=P(x)+2d/dx K(x,x) $, and $ \psi _1 (x; \lambda_k , \lambda_k ) ,
\ldots , \psi _{m_k} (x; \lambda_k , \lambda_k )$ are linearly independent.
In
addition, it also satisfies the following initial conditions: \begin{eqnarray}
       \psi _l (0; \lambda , \lambda_k ) &=&\mbox{} B^* \theta _l (k) ,
       \label{Th2.4-3} \\
       \psi _l '(0; \lambda , \lambda_k ) &=&
               \mbox{} (-A^* + K(0,0)B^* )\theta _l (k) \label{Th2.4-4} \\
         &=&\mbox{}(-A^* -B^*( \sum_{r=0}^\infty \sum_{i=1}^{m_r}
                             c_r^i \theta _i(r) \theta _i^*(r) )BB^*)
                   \theta _i (k) , \nonumber
\end{eqnarray}
or, equivalently,
\begin{equation}  \label{Th2.4-5}
          B \psi _l '(0; \lambda , \lambda_k )
           +\tilde{A} \psi _l (0; \lambda , \lambda_k )={\bf 0},
\end{equation}
where
\begin{equation} \label{Th2.4-6}
            \tilde{A} = A-BK(0,0).
\end{equation}
\end{c4}

\noindent {\bf Proof.} By (\ref{Th2.4-1}), we have
\begin{eqnarray} \nonumber
  \psi _l '' (x; \lambda , \lambda_k ) &=&\mbox{}\phi _l '' (x; \lambda ,
  \lambda_k )+[ \int _0^x K(x,t) \phi _l (t; \lambda , \lambda_k ) dt ]'' \\
   &=&\mbox{}(P(x)-\lambda I)\phi _l (x; \lambda , \lambda_k )
      +\int _0^x \frac {\partial ^2}{\partial x^2 }K(x,t)
                   \phi _l (t; \lambda , \lambda_k )dt \label{Th2.4-7}  \\
   & &\mbox{}+ K(x,x) \phi _l '(x; \lambda , \lambda_k )
             + \frac {\partial }{\partial x }K(x,t)|_{t=x}
                \phi _l (x; \lambda , \lambda_k ) , \nonumber
\end{eqnarray}
and
\begin{eqnarray*}
   \lambda \psi _l (x; \lambda , \lambda_k )
         &=& \mbox{} \lambda \phi _l (x; \lambda , \lambda_k )
             + \int _0^x K(x,t) \lambda \phi _l (t; \lambda , \lambda_k ) dt
     \\
    &=&\mbox{} \lambda \phi _l (x; \lambda , \lambda_k )
    +\int _0^x K(x,t) \lambda Y (t; \lambda )\theta _i (k) dt  \\
    &=&\mbox{}\lambda \phi _l (x; \lambda , \lambda_k )
     +\int _0^x K(x,t)  [-Y ''(t; \lambda )+P(t)Y(t;\lambda ) ]\theta _i
(k) dt  \\
    &=&\mbox{}\lambda \phi _l (x; \lambda , \lambda_k )
       +\int _0^x K(x,t) P(t) \phi _l (t; \lambda , \lambda_k ) dt \\
    & &\mbox{}   -\int _0^x K(x,t)Y''(t;\lambda )\theta _i (k)dt .
\end{eqnarray*}
Using integration by parts twice on the last integral, we have
\begin{eqnarray}
   \lambda \psi _l (x; \lambda , \lambda_k )
    &=& \lambda \phi _l (x; \lambda , \lambda_k )
       +\int _0^x K(x,t)P(t) \phi _l (t; \lambda , \lambda_k )dt
     \nonumber \\
    & &\mbox{}+(-K(x,0)A^* - \frac {\partial }{\partial y} K|_{y=0}B^* )
               \theta _l (k) \label{Th2.4-8} \\
    & &\mbox{}-K(x,x) \phi _l '(x; \lambda , \lambda_k )
        +\frac {\partial }{\partial y} K|_{y=x}
                  \phi _l (x; \lambda , \lambda_k ) \nonumber \\
    & &\mbox{}-\int _0^x \frac {\partial ^2}{\partial t^2} K(x,t)
                \phi _l (t; \lambda , \lambda_k ) dt. \nonumber
\end{eqnarray}
Then, by {\bf Theorem \ref{Th3}}, we have
\begin{eqnarray*} \lefteqn{
    -\psi _l '' (x; \lambda , \lambda_k ) +
   (\lambda I -Q(x)) \psi _l (x; \lambda , \lambda_k )\mbox{} \hskip 1cm } \\
  & &\mbox{\hskip 0.75cm}= (-K(x,0)A^* - \frac {\partial }{\partial y}
K|_{y=0}B^* )            \theta _l (k)  \\
  & &\mbox{\hskip 1.25cm} + \int _0^x
               [ \frac {\partial ^2}{\partial x^2} K
                 -\frac {\partial ^2}{\partial t^2} K -Q(x)K+KP(t)]
               dt \\
  & &\mbox{\hskip 0.75cm}={\bf 0. }
\end{eqnarray*}
Besides, for $ 1 \le l \le m_k $,
\[  \psi _l (0; \lambda , \lambda_k ) =\phi _l (0; \lambda , \lambda_k )
           =B^* \theta _l (k) , \]
\begin{eqnarray*}
 \psi _l ' (0; \lambda , \lambda_k )&=&\phi _l ' (0; \lambda , \lambda_k )
             +K(0,0)\phi _l (0; \lambda , \lambda_k )  \\
      &=&  \mbox{} (-A^* + K(0,0)B^* )\theta _l (k)   \\
      &=&\mbox{}(-A^* -B^*( \sum_{r=0}^\infty \sum_{i=1}^{m_r}
                     c_r^i \theta _i(r) \theta _i^*(r) )BB^*)\theta _l (k) .
\end{eqnarray*}
If we denote $\tilde{A} =A-BK(0,0)$, then $B \tilde{A}^*=\tilde{A}B^* $, and,
by using $BA^*=AB^*$, we have
\begin{eqnarray*} \lefteqn{
     B \psi _l ' (0; \lambda , \lambda_k ) +\tilde{A} \psi
    _l(0; \lambda , \lambda_k ) }\\
    & &\mbox{ \hskip 0.7cm}= B( -A^* + K(0,0)B^* )\theta _l (k)
        +(A-BK(0,0))B^* \theta _l (k)   \\
    & &\mbox{ \hskip 0.7cm}= {\bf 0 } .
\end{eqnarray*}
The linear independence of those vector-valued functions $\psi _l (x;
\lambda_k , \lambda_k ) , 1 \le l \le m_k $ can be proven by using (\ref{Th2.4-1}), Gronwall's
lemma and the linear independence of those functions $\phi _l (x; \lambda_k ,
\lambda_k ), 1\le l \le m_k $. $\Box $
\vskip 0.5cm
At the end of this section, we state a theorem which indicates the possible
candidates for eigenfunctions of the isospectral problem $(Q,\tilde{A},B, \tilde{\cal A}, {\cal B})$
of $(P,A,B,{\cal A,B})$ which shall be
described in next section. Furthermore, we may also use this theorem to
construct that an isospectral Dirichlet's problem of a given Dirichlet's
problem. \par

\newtheorem{c5}[c1]{Theorem} \begin{c5} \label{Th5}
Suppose $\lambda _k  \in \sigma ( P ,A ,B, {\cal A}, {\cal B} ) $, $ k \ge 0 $.
Then
\begin{eqnarray}  \label{Th2.5-1}
  \psi _l (x; \lambda_k , \lambda_k )&=&\phi _l  (x; \lambda_k , \lambda_k ) \\
      & &\mbox{} -\sum _{r=0}^\infty \sum _{i=1}^{m_r}
        c_r^i \psi _i  (x; \lambda _r , \lambda_r )
        \int _0^x \phi _i^*  (t; \lambda_r , \lambda_r )
        \phi _l (t; \lambda _k , \lambda_k )dt   \nonumber
\end{eqnarray}
for all $ 1 \le l \le m_k $.
\end{c5}

\noindent{\bf Proof.} The proof is similar to the one of [3, {\bf Theorem 1.4}].
Denote $\phi_{i,k} (x) = \phi_i (x; \lambda_k ,\lambda_k )$.
By (\ref{F}), and the integral equation (\ref{IE}), we have
\begin{eqnarray}
    K(x,y) &=& - {\cal F}(x,y) - \int _0^x K(x,t){\cal F}(t,y) dt \nonumber \\
           &=& - \sum _{r=0}^\infty \sum _{i=1}^{m_r}
                 c_r^i [\phi_{i,r} (x)+ \int_0^x K(x,t)\phi_{i,r} (t)dt]\phi_{i,r}^* (y) .
            \nonumber
\end{eqnarray}
Then, by (\ref{Th2.4-1}), we have
\begin{equation} \label{Th2.5-2}
  K(x,t)= -\sum _{r=0}^\infty \sum _{i=1}^{m_r}
            c_r^i \psi _i  (x; \lambda _r , \lambda_r )
             \phi _i^T  (t; \lambda_r , \lambda_r ).   \end{equation}
Apply (\ref{Th2.5-2}), (\ref{Th2.4-1}) implies (\ref{Th2.5-1}). $\Box $  \par

\section{ Isospectral problem }
\par
Those theorems in previous section enable us to construct an
isospectral problem from a given eigenvalue problem $(P,A,B,{\cal A,B})$ and a sequence
of real numbers $c_k^i, 1\le i \le m_k , k\ge 0 $, where the sequence $(c_k^i )$
satisfies the assumption of {\bf Theorem \ref{Th2.2}}.
As {\bf Theorem \ref{Th4}}
states, for any $k \ge 0 $, and for each $l$, $1\le l \le m_k $ , the
vector-valued
function $\psi _l (x; \lambda_k , \lambda_k )$ satisfies the boundary condition
\[     B\psi '(0)+ \tilde{A}\psi (0) = {\bf 0}, \]
where $\tilde{A}$ is given by (\ref{Th2.4-6}). Hence, the final step for
constructing isospectral problem is to determine the form of boundary condition
to be satisfied at $x=\pi $. For this purpose we use formula (\ref{Th2.5-1}).
  By (\ref{Th2.5-1}), we have
\begin{equation} \label{a3.1}
     \psi _l (\pi ; \lambda_k , \lambda_k )= \frac
{\phi _l (\pi ; \lambda_k , \lambda_k )}{1+c_k^l||\phi _l (x; \lambda_k ,
\lambda_k )||^2 }.
\end{equation}
   Differentiating (\ref{Th2.5-1}) with respect to $x$ and evaluating it at
$\pi $, we have
\begin{eqnarray*}
 \psi _l ' (\pi ; \lambda_k , \lambda_k )&=&
    \phi _l '(\pi ; \lambda_k , \lambda_k ) - c_k^l
          \psi _l '(\pi ; \lambda_k , \lambda_k ) ||\phi _l (x; \lambda_k ,
   \lambda_k )||^2 \\
   & &\mbox{}- [\sum_{r=0}^\infty \sum_{i=1}^{m_r} c_r^i
               \psi _i (\pi ; \lambda_r , \lambda_r )
               \phi _l^* (\pi ; \lambda_r , \lambda_r )]
               \phi _l (\pi ; \lambda_k , \lambda_k ) ,
\end{eqnarray*}
and, hence we have
\begin{eqnarray*}
\lefteqn{(1+c_k^l||\phi _l (x; \lambda_k ,\lambda_k )||^2 )
         \psi _l ' (\pi ; \lambda_k , \lambda_k )} \\
&=& \phi _l ' (\pi ; \lambda_k , \lambda_k )
            -[\sum_{r=0}^\infty \sum_{i=1}^{m_r}
              \frac { c_r^i \phi _i (\pi ; \lambda_r, \lambda_r )
                          \phi _i^* (\pi ; \lambda_r, \lambda_r )}{
                    1+c_r^i||\phi _i (x; \lambda_r ,\lambda_r )||^2}]
              \phi _l (\pi ; \lambda_k , \lambda_k )
\end{eqnarray*}
Acting on the above identity by ${\cal B}$ and using the condition $ {\cal B}
\phi_l' (\pi ; \lambda_k , \lambda_k ) + {\cal A}\phi_l (\pi ; \lambda_k , \lambda_k )=0 $,we have
\begin{eqnarray*}
\lefteqn{{\cal B}(1+c_k^l||\phi _l (x; \lambda_k ,\lambda_k )||^2 )
         \psi _l ' (\pi ; \lambda_k , \lambda_k )} \\
&=& {\cal B} \phi _l ' (\pi ; \lambda_k , \lambda_k )
   -{\cal B} [\sum_{r=0}^\infty \sum_{i=1}^{m_r}
              \frac { c_r^i \phi _i (\pi ; \lambda_r, \lambda_r )
                          \phi _i^* (\pi ; \lambda_r, \lambda_r )}{
                    1+c_r^i||\phi _i (x; \lambda_r ,\lambda_r )||^2}]
              \phi _l (\pi ; \lambda_k , \lambda_k )   \\
&=& -({\cal A}+{\cal B}
                [\sum_{r=0}^\infty \sum_{i=1}^{m_r}
                 \frac { c_r^i \phi _i (\pi ; \lambda_r, \lambda_r )
                             \phi _i^* (\pi ; \lambda_r, \lambda_r )}{
                       1+c_r^i||\phi _i (x; \lambda_r ,\lambda_r )||^2}] )
               \phi _l (\pi ; \lambda_k , \lambda_k ).
\end{eqnarray*}
Then, by (\ref{a3.1}), we have
\[
   {\cal B} \psi _l ' (\pi ; \lambda_k , \lambda_k )
          = -({\cal A+B}
                 [\sum_{r=0}^\infty \sum_{i=1}^{m_r}
                  \frac { c_r^i \phi _i (\pi ; \lambda_r, \lambda_r )
                              \phi _i^* (\pi ; \lambda_r, \lambda_r )}{
                        1+c_r^i||\phi _i (x; \lambda_r ,\lambda_r )||^2}] )
              \psi _l (\pi ; \lambda_k , \lambda_k ) ,
\]
and hence,
\[
   {\cal B} \psi _l ' (\pi ; \lambda_k , \lambda_k )
  +\tilde{{\cal A}} \psi _l  (\pi ; \lambda_k , \lambda_k ) = {\bf 0} ,, \]
where
\begin{equation}  \label{a3.2}
     \tilde{{\cal A}}= {\cal A+B}
               [\sum_{r=0}^\infty \sum_{i=1}^{m_r}
                \frac { c_r^i \phi _i (\pi ; \lambda_r, \lambda_r )
                            \phi _i^* (\pi ; \lambda_r, \lambda_r )}{
                      1+c_r^i||\phi _i (x; \lambda_r ,\lambda_r )||^2}] .
\end{equation}
In fact, by (\ref{Th2.5-2}), (\ref{a3.2}) can be simplified as
\begin{equation} \label{a3.3}
      \tilde{\cal A} = {\cal A}- {\cal B} K(\pi ,\pi) .
\end{equation}
Furthermore, if we can prove that, for any ${\bf R}^N$-valued function $f$
satisfying $ Bf'(0)+ \tilde{A} f(0) ={\bf 0} $ , $ {\cal B} f' (\pi )
+ \tilde{\cal A} f (\pi ) ={\bf 0} $, and
\[
      \int_0^\pi \langle f(x), \psi _l ( x ; \lambda _k , \lambda _k ) \rangle
dx        =0, \hskip 0.5cm k \ge 0, \hskip 0.25cm l=1, \ldots , m_k,
\] we have $f \equiv {\bf 0}$, then the set $ \{ \psi _l ( x ; \lambda _k ,
\lambda _k ): k\ge 0, 1\le l \le m_k \} $ is complete, and, by {\bf Theorem
\ref{Th4}}, we have $\Sigma (P, A, B, {\cal A, B}) = \Sigma (Q, \tilde{A} , B,
\tilde{\cal A} , {\cal B})$. For our purpose, let
\[       T(f)(x)= \int _0^x
                      K(x,t) f(t) dt. \]
Then
\[
      T^* (g)(t)= \int _t^\pi K^*(x,t)g(x)dx.
\]
Writing $ \psi _l ( x ; \lambda _k , \lambda _k ) = (I+T)
               \phi _l ( x ; \lambda _k , \lambda _k ) $ , $ 1\le l \le m_k $,
we have
\begin{eqnarray*}
  0= \int_0^\pi \langle f(x), \psi _l ( x ; \lambda _k , \lambda _k ) \rangle
dx   &=&\int_0^\pi \langle f(x), (I+T)
          \phi _l ( x ; \lambda _k , \lambda _k ) \rangle dx  \nonumber \\
  &=&\int_0^\pi \langle (I+T^*)f(x),
                       \phi _l ( x ; \lambda _k , \lambda _k ) \rangle dx.
\end{eqnarray*}
Now set $g= (I+T^*)f $. If we show that $ Bg'(0)+ Ag(0) ={\bf 0}$, and ${\cal
B}g'(\pi) +{\cal A}g(\pi )={\bf 0} $, then by the completeness of the set $\{
\phi _l ( x ; \lambda _k , \lambda _k ) : k \ge 0, 1\le l \le m_k \} $, we
have $g \equiv {\bf 0} $ and hence $f \equiv {\bf 0} $. We only check
the identity $Bg'(0)
+Ag(0)={\bf 0} $, the other part can be proved by similar argument. Suppose
$Bf'(0)+\tilde{A} f(0) = {\bf 0} $. Then, using (\ref{Th3-2}) and
(\ref{Th2.4-6}), we have
\begin{eqnarray*}
    Bg'(0)+Ag(0) &=& B[ f'(0)-K^*(0,0)f(0)   \\
          & &\mbox{}  + \int _0^\pi K^*_t (x,0) f(x)dx ]
        + A[f(0)+\int _0^\pi K^* (x,0) f(x) dx ] \\
     &=& Bf'(0)-BK^* (0,0)f(0)  \\
      & &\mbox{}+Af(0)+ \int _0^\pi [BK^*_t (x,0) +AK^* (x,0) ]f(x) dx \\
     &=& Bf'(0)+ [A-BK(0,0)]f(0)= Bf'(0)+\tilde{A} f(0)= {\bf 0 }.
\end{eqnarray*}

As a conclusion of the previous arguments, we have the following theorem.
\newtheorem{c6}{Theorem}[section] \begin{c6} \label{Th6} Let $m_k$ denote the
multiplicity of $\lambda _k $ in $\sigma ( P, A, B, {\cal A, B} )$.
Suppose  $\{ c_k^i , 1\le i \le m_k , k\ge 0 \} $ is a sequence,
satisfying the condition (\ref{Th2-c}) and making ${\cal F}(x,y)$ in (\ref{F})
a $C^2$-function, $Q(x)$ is as that defined in {\bf Theorem 2.3}, $ \tilde{A} $ and $\tilde{
\cal A} $ are as those defined in (\ref{Th2.4-6}) and (\ref{a3.3}). Then $\Sigma
(P,
A, B, {\cal A , B} ) = \Sigma (Q, \tilde{A}, B, \tilde{\cal A} ,{\cal B}) $.
\end{c6}

As a final remark, we note that if $A=I, B=0, {\cal A}=I$, and ${\cal B}=0 $ in
(\ref{a1.1}), then the matrices $\tilde{A}$ and $\tilde{\cal A}$ in
{\bf Theorem \ref{Th6} } are equal to $I$, the identity matrix. Hence,
for a given Dirichlet problem $(P, I,0,I,0)$, the isospectral problem constructed in {\bf
Theorem \ref{Th6}} is also a Dirichlet problem.

\section{ Examples }

In this section, we use our theory to construct some examples which
have some significant meaning the scalar case can not tell.

Suppose $\lambda_\circ $ is an eigenvalue of (\ref{a1.1})
with multiplicity $m_\circ $. Let $ \phi_\circ (x) = \mbox{col}
( \phi _1 (x), $ \newline $ \phi_2 (x), \cdots, \phi_N (x)) $ 
be an eigenfunction corresponding to $\lambda _\circ $.
Take  
\begin{eqnarray*}
     {\cal F} (x,y) &= & c \phi _\circ (x) \phi_\circ^* (y)   \\
                    &= & c \left ( \begin{array}{cccc}
      \phi_1 (x) \phi_1 (y) &\phi_1 (x) \phi _2(y)  & \cdots & \phi_1 (x) \phi_N (y) \\
      \phi_2 (x) \phi_1 (y) &\phi_2 (x) \phi _2(y)  & \cdots & \phi_2 (x) \phi_N (y) \\
      \vdots & \vdots & \ddots &\vdots  \\
      \phi_N (x) \phi_1 (y) &\phi_N (x) \phi _2(y)  & \cdots & \phi_N (x) \phi_N (y) \\
      \end{array}  \right ) .
\end{eqnarray*}

Plugging it into (\ref{IE}), and letting $k_{ij} (x,y) $ denote the $(i,j)$ entry of $K(x,y)$, we have
\begin{equation}
    k_{ij} (x,y) + c \phi_i (x) \phi_j (y)
               +(c \int _0^x (\sum _{r=1}^N k_{ir} (x,t) \phi _r (t))dt)
                \phi _j (y) =0 ,
\end{equation}
 for $ i=1,\ldots ,N$ and $j=1,\ldots ,N $.

We shall show that
\begin{equation}\label{a4.1}
    k_{ij} (x,y) =- \frac {c \phi_i (x) \phi_j (y)}{ 1+ c \int _0^x |\phi _\circ (t)|^2 dt } .
\end{equation}
 For $i$ fixed, consider the equations
\begin{equation} \label{a4.2}
    k_{ij} (x,y) + c \phi_i (x) \phi_j (y)
               +(c \int _0^x (\sum _{r=1}^N k_{ir} (x,t) \phi _r (t))dt)
                \phi _j (y) =0,
\end{equation}
$ 1\le j \le N$.
Multiplying the $j$th equation by $\phi _j (y) $, integrating it from $0$
to $x$ with respect to $y$,
and denoting
$ p_{i,j} (x) = \int _0^x k_{ij} (x,t) \phi _j (t) dt $ and $ \alpha _j (x) =
\int _0^x \phi_j^2 (t)dt $, $ 1\le j \le N $,
we have the following lineaar system of equations with unknowns $p_{ij} (x) $,
\begin{equation}
  p_{ij} (x) + c \phi _i (x) \alpha _j (x) + c \alpha _j (x)( \sum _{r=1}^N
    p_{ir} (x) ) =0, \quad 1 \le j \le N .
\end{equation}
Solving this system, we obtain
\begin{equation} \label{a.4.4}
   p_{ij} (x)= - \frac {c \phi _i (x) \alpha _j (x) }{1+ c \int _0^x |\phi _\circ (t)|^2 dt }, \quad 1 \le j \le N.
\end{equation}
Note that, by (\ref{a4.2}),
\[
    k_{ij} (x,y) =- c \phi_i (x) \phi_j (y)
               -(c \int _0^x (\sum _{r=1}^N k_{ir} (x,t) \phi _r (t))dt)
                \phi _j (y) .
\]
Hence if we plug (\ref{a.4.4}) into (\ref{a4.2}), we obtain (\ref{a4.1})
.
It also follows from (\ref{a4.1}), that
\[
    K(x,y) = -\frac 1{1+ c \int _0^x |\phi _\circ (t)|^2 dt }{\cal F}(x,y).
\]
Hence, according to {\bf Theorem \ref{Th6}}, by setting
\begin{eqnarray}
        Q(x) &= & P(x) - 2 \frac d{dx} [ \frac {\phi _\circ (x) \phi_\circ^* (x)}
        {1+ c \int _0^x |\phi _\circ (t)|^2 dt } ], \nonumber \label{a4.3} \\
        \tilde{A} &=& A-cB\phi _\circ (0) \phi_\circ^* (0), \\
        \tilde{\cal A}&=& {\cal A} -\frac {c\phi _\circ (\pi) \phi_\circ^* (\pi)}{1+ c ||\phi _\circ ||^2}, \nonumber
\end{eqnarray}
we have
$\Sigma (P,A, B, {\cal A , B} ) =
 \Sigma (Q, \tilde{A}, B, \tilde{\cal A} ,{\cal B}) $.

As an example to the above construction, we construct the following  eigenvalue
problem which has an eigenvalue of multiplicity $2$.

 Let $I$ be the $2 \times 2$ identity matrix. Take
\[
   P(x) = \left ( \begin{array}{cc}
                  -3 & 0 \\
                   0 & 0
                  \end{array} \right ) , \mbox{\hskip 0.25cm} A={\cal A}= I,\mbox{\hskip 0.25cm} B={\cal B}=0.
\]
Then one can verify that for the eigenvalue problem $ (P,I,0,I,0)$,
$1$ is an eigenvalue of multiplicity $2$, and
 the other eigenvalues are all simple. The eigenspace corresponding to $1$ is
the vector space spanned by the two vector-valued functions $ ( \sin (2x),0 )^*$ and $(0, \sin (x) )^* $.
Choosing $ (\sin(2x) ,\sin(x))^*$ as the eigenfunction which plays the role of $\phi_\circ (x)$ in the above construction,
$c=1$, and using (\ref{a4.3}),
we have
\begin{eqnarray*}
    Q(x)= \left ( \begin{array}{cc}
                  -3 & 0 \\
                   0 & 0
                  \end{array} \right )&-& \frac d{dx}[ \frac {2}{1+  \int _0^x (\sin^2 (t) + sin^2 (2t))dt } \\
                    & &\cdot  \left ( \begin{array}{cc}
                                  \sin^2 (2x) & \sin (2x) \sin (x) \\
                                  \sin (x) \sin (2x) & \sin ^2 (x)
                                 \end{array}
                         \right ) ] .
\end{eqnarray*}
and $ \tilde{A} = \tilde{\cal A}=I $. Note that the matrix potential function
$Q(x)$ is not simultaneously diagonalizable since, as checked by computation,
the matrix-valued functions $Q(x)$ and $Q' (x) $ do not commute.
On the other hand, if we take $(0, \sin (x))^* $ instead, then we find that
$ Q (x) $ is a diagonal matrix-valued function and is of the following form
\[
    Q(x) = \left ( \begin{array}{cc}
                  -3 & 0  \\ \mbox{} & \mbox{} \\
                   0 &  \frac d{dx} (\frac {-2 \sin ^2 (x)}{ 1+  \int _0^x \sin^2 (t) dt } )
                 \end{array} \right ).
\]

There are lots of interesting phenomena can be observed from our construction,
which shall be observed later.

\noindent {\bf Acknowledgements. } (i) The author show his gratitude to his
Ph D. adviser Professor C. L. Shen for his instruction. (ii) 
The author became aware that Professor B. M.
Levitan and Max. Jodeit also obtained analogous results in \cite{JL3}.

\end{document}